\documentclass[reqno]{amsart}
\usepackage{amssymb,url,verbatim}
\date{}
\allowdisplaybreaks[4]
\theoremstyle{plain}
\newtheorem{thm}{Theorem}
\newtheorem{lem}{Lemma}

\theoremstyle{remark}

\theoremstyle{definition}

\DeclareMathOperator{\td}{d\!}
\newcommand{\tn}{\mathbb{N}}

\newcommand{\nbd}{\nobreakdash}
\usepackage[colorlinks,linkcolor=blue,urlcolor=red,linktocpage]{hyperref}
\begin{document}
\title
{The best bounds of harmonic sequence}

\author[Ch.-P. Chen]{Chao-Ping Chen}
\address[Ch.-P. Chen]{Department of Applied Mathematics and Informatics, Jiaozuo Institute of Technology, Jiaozuo City, Henan 454000, CHINA}

\author[F. Qi]{Feng Qi}
\address[F. Qi]{Department of Applied Mathematics and Informatics, Jiaozuo Institute of Technology, Jiaozuo City, Henan 454000, CHINA}
\email{qifeng@jzit.edu.cn, fengqi618@member.ams.org}
\urladdr{http://rgmia.vu.edu.au/qi.html}

\begin{abstract}
For any natural number $n\in\tn$,
\begin{equation}
\frac{1}{2n+\frac1{1-\gamma}-2}\le \sum_{i=1}^n\frac1i-\ln n-\gamma<\frac{1}{2n+\frac13},
\end{equation}
where $\gamma=0.57721566490153286\dotsm$ denotes Euler's constant. The constants $\frac{1}{1-\gamma}-2$ and $\frac13$ are the best possible. 
\par
As by-products, two double inequalities of the digamma and trigamma functions are established.
\end{abstract}

\keywords{Inequality, best bound, harmonic sequence, Binet formula, digamma function, trigamma function}

\subjclass[2000]{26D15, 33B15}

\thanks{The authors were supported in part by NNSF (\#10001016) of China, SF for the Prominent Youth of Henan Province (\#0112000200), SF of Henan Innovation Talents at Universities, Doctor Fund of Jiaozuo Institute of Technology, CHINA}

\thanks{This paper was typeset using \AmS-\LaTeX}

\maketitle	

\section{Introduction}

Let $n$ be a natural number, then we have
\begin{equation}
\label{franel}
\frac1{2n}-\frac1{8n^2}<\sum_{i=1}^n\frac1i-\ln n-\gamma<\frac1{2n},
\end{equation}
where $\gamma=0.57721566\dotsm$ is Euler's constant.
\par
The inequality \eqref{franel} is called in literature Franel's inequality \cite[Ex.~18]{PS}. Because of the well known importance of the harmonic sequence $\sum_{i=1}^n\frac1i$, there exists a very rich literature on inequalities of the harmonic sequence $\sum_{i=1}^n\frac1i$. For example, \cite{abram,harijmest}, \cite[pp.~68\nbd--78]{kuang} and references therein.
\par
L. T\'oth and S. Mare in \cite[p. 264]{laszl} proposed the following problems:
\begin{enumerate}
\item 
Prove that for every positive integer $n$ we have
\begin{equation}
\label{laszl}
\frac{1}{2n+\frac25}<1+\frac12+\cdots +\frac{1}{n}-\ln n-\gamma<\frac{1}{2n+\frac13},
\end{equation}
where $\gamma$ is Euler's constant.
\item
Show that $\frac25$ can be replaced by a slightly smaller number, but that $\frac13$ cannot be replaced by a slightly larger number.
\end{enumerate}
\par
In \cite{yanghar}, basing on improving of Euler-Maclaurin summation formula, some general inequalities of the harmonic sequence $\sum_{i=1}^n\frac1i$ are established, including recovery of inequality \eqref{laszl}. In 1999, Sh.\nbd-R.~Wei and B.\nbd-Ch.~Yang in \cite{wei} verified inequality \eqref{laszl} again by calculus.
\par
In 1997, K.~Wu and B.\nbd-Ch.~Yang in \cite{wuk} proved the second problem due to T\'oth and Mare by using the following
\begin{equation}
\label{approx}
\sum_{i=1}^n\frac1i=\gamma+\ln n+\frac1{2n}-\sum_{i=1}^{q-1}\frac{B_{2i}}{2in^{2i}}-\int_n^\infty \frac{B_{2q}(x)}{x^{2q}}\td x
\end{equation}
and
\begin{equation}
\label{berhar}
\int_n^\infty \frac{B_{2q-1}(x)}{x^{2q}}\td x<\frac{(-1)^qB_{2q}}{2qn^{2q}},
\end{equation}
where $B_i(x)$ is Bernoulli's polynomial, $B_{2i}=B_{2i}(n)$ is Bernoulli's number for $i\in\tn$, and $n$ and $q$ are positive integers,.
\par
In this short note, through establishing, by exploiting the well known first Binet's formula, two double inequalities of the digamma and trigamma functions and by utilizing the approximating expansion of digamma function $\psi$, the best lower and upper bounds of the sequence $\sum_{i=1}^n\frac1i-\ln n-\gamma$ are given. Thses results refine inequality \eqref{laszl} and solve affirmatively the second problem mentioned above.

\begin{thm}\label{franbest}
For any natural number $n\in\tn$, we have
\begin{equation}\label{ding}
\frac{1}{2n+\frac1{1-\gamma}-2}\le \sum_{i=1}^n\frac1i-\ln n-\gamma<\frac{1}{2n+\frac13},
\end{equation}
where $\gamma=0.57721566490153286\dotsm$ denotes Euler's constant. The constants $\frac{1}{1-\gamma}-2$ and $\frac13$ are the best possible. 
\end{thm}

\section{Lemma}

In order to prove inequality \eqref{laszl}, the following lemma is necessary.

\begin{lem}
For $x>0$, we have
\begin{gather}\label{lem1}
\frac{1}{2x}-\frac{1}{12x^2}<\psi(x+1)-\ln x<\frac{1}{2x}\\
\intertext{and}
\label{lem2}
\frac{1}{2x^2}-\frac{1}{6x^3}<\frac{1}{x}-\psi '(x+1)<\frac{1}{2x^2}-\frac{1}{6x^3}+\frac{1}{30x^5}, 
\end{gather}
where $\psi= \frac{\Gamma'}{\Gamma}$ is the logarithmic derivative of the gamma function 
\begin{equation} 
\Gamma(x)=\int_0^\infty t^{x-1}e^{-t}\td t.
\end{equation}
\end{lem}

\begin{proof}
It is a well known fact (\cite{abram} and \cite[p.~103]{wang}) that for $x>0$ and a nonnegative integer $m$, 
\begin{equation} \label{gamm}
\psi(x+1)=\psi(x)+\frac{1}{x}
\end{equation}
and
\begin{equation}
\label{gam}
\frac{m!}{x^{m+1}}=\int_{0}^{\infty}t^{m}e^{-xt}\td t.
\end{equation}
\par
The first Binet's formula (\cite{abram} and \cite[p.~106]{wang}) states that for $x>0$
\begin{equation}\label{binet}
\ln \Gamma(x)= \left(x-\frac{1}{2}\right)\ln x-x+\ln \sqrt{2\pi}-\int_{0}^{\infty}\left(\frac{1}{2}+\frac{1}{t}-\frac{1}{1-e^{-t}}\right)\frac{e^{-xt}}{t}\td t.
\end{equation}
Differentiating \eqref{binet}, integrating by part and using formulas \eqref{gam} and \eqref{gamm}, it is deduced that
\begin{equation}
\label{binetcor}
\psi(x+1)-\ln x
=\int_{0}^{\infty}\left(\frac{1}{t}-\frac{1}{e^t-1}\right)e^{-xt}\td t.
\end{equation}
Using formulas \eqref{gam} and \eqref{binetcor} and the series expansion of $e^x$ at $x=0$ yields
\begin{equation}
\begin{split}
& \quad \psi(x+1)-\ln x-\frac{1}{2x}+\frac{1}{12x^2} \\
&=\int_{0}^{\infty}\left(\frac{1}{t}-\frac{1}{e^t-1}-\frac{1}{2}+\frac{1}{12}t \right)e^{-xt}\td t \\
&=\int_{0}^{\infty}\frac{12(e^t-1)-12t-6t(e^t-1)+t^2(e^t-1)}{12t(e^t-1)}e^{-xt}\td t \\
&=\int_{0}^{\infty}\biggl[\frac{1}{12t(e^t-1)}\sum_{n=5}^{\infty}\frac{(n-3)(n-4)}{n!}t^n \biggr]e^{-xt}\td t\\
&>0
\end{split}
\end{equation}
and
\begin{equation}
\begin{split}
\psi(x+1)-\ln x-\frac{1}{2x}
&=\int_{0}^{\infty}\left(\frac{1}{t}-\frac{1}{e^t-1}-\frac{1}{2} \right)e^{-xt}\td t \\
&= -\int_{0}^{\infty}\biggl[\frac{1}{2t(e^t-1)}\sum_{n=3}^{\infty}\frac{n-2}{n!}t^n \biggr]e^{-xt}\td t \\
&<0.
\end{split}
\end{equation}
Hence, inequality \eqref{lem1} follows.
\par
Differentiation of \eqref{binetcor} immediately produces
\begin{equation}
\label{binetcorcor}
\frac{1}{x}-\psi '(x+1)=\int_{0}^{\infty}\left(1-\frac{t}{e^t-1}\right)e^{-xt}\td t.
\end{equation}
Exploiting formulas \eqref{gam} and  \eqref{binetcorcor} and the series expansion of $e^x$ at $x=0$ yields
\begin{equation}
\begin{split}
& \quad \frac{1}{x}-\psi '(x+1)-\frac{1}{2x^2}+\frac{1}{6x^3} \\
&=\int_{0}^{\infty}\left(1-\frac{t}{e^t-1}-\frac{1}{2}t+\frac{1}{12}t^2 \right)e^{-xt}\td t \\
&=\int_{0}^{\infty}\biggl[\frac{1}{12(e^t-1)}\sum_{n=5}^{\infty}\frac{(n-3)(n-4)}{n!}t^n \biggr]e^{-xt}\td t\\
&>0
\end{split}
\end{equation}
and
\begin{equation}
\begin{split}
& \quad \frac{1}{x}-\psi '(x+1)-\frac{1}{2x^2}+\frac{1}{6x^3}-\frac{1}{30x^5} \\
&=\int_{0}^{\infty}\left(1-\frac{t}{e^t-1}-\frac{1}{2}t+\frac{1}{12}t^2-\frac{1}{720}t^4 \right)e^{-xt}\td t \\
&= \int_{0}^{\infty}\biggl[\frac{1}{720(e^t-1)}\sum_{n=7}^{\infty}\left(\frac{720}{n!}-\frac{360}{(n-1)!}+\frac{60}{(n-2)!}-\frac{1}{(n-4)!}\right)t^n \biggr]e^{-xt}\td t.
\end{split}
\end{equation}
Noticing that for $n\ge7$, 
\begin{equation}
\begin{split}
& \quad  \frac{720}{n!}-\frac{360}{(n-1)!}+\frac{60}{(n-2)!}-\frac{1}{(n-4)!}	  \\
&= -\frac{120+218(n-7)+119(n-7)^2+22(n-7)^3+(n-7)^4}{n!}<0,
\end{split}
\end{equation}
we obtain
\begin{equation}
 \frac{1}{x}-\psi '(x+1)-\frac{1}{2x^2}+\frac{1}{6x^3}-\frac{1}{30x^5}<0.
\end{equation}
Therefore, inequality \eqref{lem2} holds. The proof is complete.
\end{proof}

\section{Proof of Theorem \ref{franbest}}

In \cite{abram}, \cite[p.~593]{shouce} and \cite[p.~104]{wang} it is given that $\psi(n)=\sum_{k=1}^{n-1}\frac{1}{k}-\gamma$. Thus, inequality \eqref{ding} can be rearranged as
\begin{equation}
\frac{1}{3}<\frac{1}{\psi(n+1)-\ln n}-2n \le \frac{1}{1-\gamma}-2.
\end{equation}
\par
Define for $x>0$
\begin{equation}
\phi(x)=\frac{1}{\psi(x+1)-\ln x}-2x.
\end{equation}
Differentiating $\phi$ and utilizing \eqref{lem1} and \eqref{lem2} reveals that for $x>\frac{12}5$,
\begin{equation}
\begin{split}
& \quad (\psi(x+1)-\ln x)^2\phi '(x)\\
&=\frac{1}{x}-\psi '(x+1)-2(\psi(x+1)-\ln x)^2 \\
&<\frac{1}{2x^2}-\frac{1}{6x^3}+\frac{1}{30x^5}-2\left(\frac{1}{2x}-\frac{1}{12x^2}\right)^2 \\
&= \frac{12-5x}{360x^5}<0,
\end{split}
\end{equation}
and $\phi(x)$ decreases with $x>\frac{12}5$.
\par
Straightforward calculation produces
\begin{align}
\phi(1)&=\frac1{1-\gamma}-2=0.36527211862544155\dotsm,\\
\phi(2)&=\frac{1}{\frac32-\gamma-\ln2}-4=0.35469600731465752\dotsm,\\
\phi(3)&=\frac{1}{\frac{11}6-\gamma-\ln3}-6=0.34898948531361115\dotsm.
\end{align}
Therefore, the sequence
\begin{equation}
\phi(n)=\frac{1}{\psi(n+1)-\ln n}-2n,  \quad n\in\tn
\end{equation}
is decreasing strictly, and for $n\in\tn$
\begin{equation}
\lim_{n \to \infty}\phi(n)<\phi(n)\le \phi(1)=\frac{1}{1-\gamma}-2.
\end{equation}
\par
Making use of approximating expansion of $\psi$ in \cite{abram}, \cite[p.~594]{shouce}, or \cite[p.~108]{wang} gives
\begin{equation}
\psi(x)=\ln x-\frac{1}{2x}-\frac{1}{12x^2}+O(x^{-4}) \quad (x \to \infty),
\end{equation}
and then
\begin{equation}											
\lim_{n \to \infty}\phi(n)=\lim_{x \to \infty}\phi(x) 
=\lim_{x \to \infty}\frac{\frac13+O(x^{-2})}{1+O(x^{-1})}=\frac{1}{3}. 
\end{equation}
The proof is complete.

\end{document}